\documentclass[a4paper]{amsart}

\usepackage{amssymb,latexsym}
\usepackage{amsmath,amscd}
\usepackage[all]{xy}

\title{\bf Prequantization and Lie brackets}
\author{Marius Crainic}
\thanks{{\tt research supported by a KNAW fellowship (Utrecht) and the Miller Foundation (Berkeley)}}
\address{Department of Mathematics\\
        Utrecht University, P.O. Box 80.010, 3508 TA\\
        Utrecht, The Netherlands\\
        crainic@math.uu.nl}
\date{}

\newcommand{\rmap}{\longrightarrow}

\newcommand{\SP} [1]     {{\left\langle {{#1}} \right\rangle}}

\newtheorem{lemma} {Lemma} [section]

\newtheorem{theorem} [lemma] {Theorem}
\newtheorem{corollary} [lemma] {Corollary}
\newtheorem{definition}[lemma] {Definition}

\newtheorem{remark}[lemma]{Remark}

\numberwithin{equation}{section}

\begin{document}

\maketitle

\begin{abstract} 
We start by describing the relationship between the classical 
prequantization condition and the integrability of a certain Lie 
algebroid associated to the problem and use this to give a global 
construction of the prequantizing bundle in terms of path spaces
(Introduction), then we rephrase the problem in terms of groupoids 
(second section), and then we study the more general problem of 
prequantizing groupoids endowed with a multiplicative 2-forms 
(the rest of the paper). 
\end{abstract}


\tableofcontents

\section{Introduction: integrality and integrability (classical prequantization)}

We start by describing one of the simplest (and motivating) particular cases
of the problem studied in this paper (classical prequantization), then we give an outline
of the main results/techniques and we briefly describe the connections with existing literature.  
The starting point for this paper is the striking similarity
between the following two known results, similarity that begs for more explanations.

\begin{theorem}
\label{theorem1}
 A simply connected\footnote{this assumption is made only to simplify the presentation} symplectic manifold $(M, \omega)$ is prequantizable if and only if $Per(\omega)= k \mathbb{Z}$, for some integer $k$. Moreover, the prequantization is unique (up to isomorphism). 
\end{theorem}

\begin{theorem}
\label{theo2} Given a closed 2-form $\omega\in \Omega^2(M)$, the Lie algebroid $A_{\omega}$ associated to $\omega$ is integrable if and only if $Per(\omega)= a \mathbb{Z}$ for some real number $a$.
\end{theorem}

Recall that the period group of a closed 2-form $\omega$, $Per(\omega)$, is the subgroup of $\mathbb{R}$ consisting of all integrals
\[ \int_{\gamma} \omega := \int_{0}^{1}\gamma^*\omega,\]
with $\gamma: S^2\rmap M$. Also, a prequantization of $\omega$ is a complex line bundle $\mathbb{L}$ together with a connection $\nabla$, so that $\omega$ coincides with the first Chern-form $c_{1}(\mathbb{L}, \nabla)$. Equivalently, this is the same as a principal $S^1$-bundle $\pi: P\rmap M$ together with a connection 1-form $\theta\in \Omega^1(P)$  such that
$\pi^{*}\omega= d \theta$. Of course, this does not require $\omega$ to be non-degenerated, and Theorem \ref{theorem1} is true
for all closed 2-forms. The statement of Theorem \ref{theo2} will be explained below. 

The main question we have in mind is to {\it find the precise relation between the prequantization of $\omega$ and the integrability of $A_{\omega}$}, and to see what integrability tells us about prequantization. We start here by describing this relation, and we 
leave the proofs for the forthcoming sections (in a much more general context).

\subsection{From prequantization to integrability} Given $P$, we consider the gauge group of $P$, $Gauge(P)$, consisting of all automorphisms of the principal bundle $P$ (with base map an arbitrary diffeomorphism of $M$). Its Lie algebra is
\[ \mathfrak{g}_{\omega}= \mathcal{X}(M)\oplus C^{\infty}(M) ,\]
with the bracket
\[ [(X, f), (Y, g)]= ([X, Y], X(f)- Y(g)+ \omega(X, Y)) .\]
As emphasized by the notation, {\it $\mathfrak{g}_{\omega}$ depends only on $\omega$, and not on the bundle $P$}. Hence, 
viewing $\omega$ as the starting data, we can re-interpret the role of $P$ as follows: {\it a prequantizing bundle $P$ of $\omega$ gives
us a Lie group (namely $Gauge(P)$) integrating  $\mathfrak{g}_{\omega}$.}

There is an even better way to state this remark, which avoids working with infinite dimensional Lie groups and Lie algebras,
and which exploits the full structure underlying our discussion. For this, we move our discussion into the world of Lie algebroids, which can be thought of as infinite dimensional Lie algebras ``of geometric type'', i.e. defined on the space of section of a vector bundle, and obeying a Leibniz-type identity. More formally, a Lie algebroid over $M$ consists of a vector bundle $A$ over $M$, together with a Lie algebra structure $[\cdot, \cdot]$ on the space of sections $\Gamma(A)$ and a bundle map $\rho: A\rmap TM$, satisfying
\[ [\alpha, f\beta]= f[\alpha, \beta]+ \mathcal{L}_{\rho(\alpha)}(f)\beta,\]
for all $\alpha, \beta\in \Gamma(A)$, $f\in C^{\infty}(M)$, where $\mathcal{L}$ stands for the Lie derivative along vector fields.
Back to our 2-form $\omega$, one has $\mathfrak{g}_{\omega}= \Gamma(A_{\omega})$, where
\[ A_{\omega}:= TM\oplus \mathbb{R}\]
is the Lie algebroid with $\rho$ being the projection on the first
component, and the bracket of $\mathfrak{g}_{\omega}$. In general,
$A_{\omega}$ defines an extension of $TM$ by the trivial line
bundle $\mathbb{R}$ (with zero anchor and zero bracket)
\[ 0\rmap \mathbb{R}\rmap A_{\omega}\rmap TM\rmap 0 ,\]
and, exactly as in the case of extensions of Lie algebras, such
extensions are classified by closed
2-forms.

As in the case of infinite dimensional Lie algebras, not all algebroids are integrable (in the sense that they come from a Lie groupoid), and the algebroids of type type $A_{\omega}$ produce the simplest non-integrable examples.  However, when prequantizing bundles exist, $A_{\omega}$ must be integrable: a prequantizing bundle induces a ``gauge groupoid'' $\mathcal{G}auge(P)$ over $M$ (which is the quotient of $P\times P$ by the action of $S^1$, see also the next section), and then, the ``oid''-version of the remark above says that {\it a prequantizing bundle $P$ of $\omega$ gives us a Lie groupoid integrating  $A_{\omega}$ (namely $\mathcal{G}auge(P)$.}

\subsection{The converse}  The conclusion above clearly explains how prequantization affects integrability.
However, it is the reverse direction that is more interesting (and probably less known). Here we explain how the general constructions related 
to the integrability of Lie algebroids translate into classical terms and give us a global construction of prequantization.
So, let $\omega$ be a closed 2-form on $M$, and we consider
\[ \mathcal{S}(\omega):= \mathbb{R}/ Per(\omega),\]
called the structural group of $\omega$.

The Banach manifold $P(M)$ of all $C^2$-paths in $M$ has a 1-form $\sigma$ naturally made out of $\omega$: given $\gamma\in P(M)$, a
tangent vector $X_{\gamma}\in T_{\gamma}P(M)$ is a path $X_{\gamma}: I\rmap TM$ over $\gamma$ (i.e. $X_{\gamma}(t)\in T_{\gamma(t)}M$) and
\[ \sigma(X_{\gamma})= \int_{0}^{1} \omega(\frac{d\gamma}{dt}, X_{\gamma}(t)) dt .\]
We will see that this is just a pull-back of the canonical contact (Liouville) form on cotangent bundles.

Next, we fix a point $x_0\in M$, and we consider the submanifold $P_{x_0}(M)$ of those paths that start at $x_0$.
On the manifold $P_{x_0}(M)\times \mathbb{R}$ we consider:
\begin{itemize}
\item  The difference between the pull-back of $\sigma$ and the pull-back of the canonical 1-form $dt$ on $\mathbb{R}$ (by the projections):
\[\tilde{\theta}:= pr_{2}^{*}(dt)- pr_{1}^{*}(\sigma) \ \ (\text{a\ $1$-\ form\ on}\ P_{x_0}(M)\times \mathbb{R}).\]
\item the equivalence relation $\sim$ given by:
\[ (\gamma_0, r_0)\sim (\gamma_1, r_1) \Longleftrightarrow r_0- r_1= \int_{\gamma} \omega \]
for some homotopy $\gamma$ between $\gamma_0$ and $\gamma_1$ with fixed end points.
\end{itemize}
We also consider the quotient space
\[ P_{\omega}:= P_{x_0}(M)\times \mathbb{R}/\sim ,\]
and we look for smooth structures on $P_{\omega}$ for which the canonical projection from $P_{x_0}(M)\times \mathbb{R}$ is a submersion.
Such a smooth structure will be unique if it exists, and in that case we say that $P_{\omega}$ is smooth. Note also that,
by acting on the second component, we obtain a free action of the structural group $\mathcal{S}(\omega)$ on $P_{\omega}$. In particular, this shows
that the smoothness of $P_{\omega}$ forces $Per(\omega)$ to be closed in $\mathbb{R}$, i.e. of type $a\mathbb{Z}$ for some $a\in \mathbb{R}$.
Even better, we have the following which 
gives a global construction of prequantizing bundles and their connection forms.

\begin{theorem}
Given a closed two-form $\omega$ on a simply connected manifold $M$, we have:
\begin{enumerate}
\item[(i)] the orbits of $\sim$ form the leaves of a regular foliation on $P_{x_0}(M)\times\mathbb{R}$, of codimension $\text{dim}(M)+ 1$.
\item[(ii)] $P_{\omega}$ is smooth if and only if $Per(\omega)= a\mathbb{Z}$ for some $a\in \mathbb{R}$. Furthermore, in this case $P_{\omega}$
is a principal $\mathcal{S}(\omega)$-bundle over $M$, and $\tilde{\theta}$ descends to a 1-form $\theta$ on $P_{\omega}$.
\item[(iii)] If $Per(\omega)= \mathbb{Z}$, then $(P, \theta)$ is the (unique) prequantization of $\omega$.
\end{enumerate}
\end{theorem}

Note that, when $Per(\omega)= k\mathbb{Z}$ for some $k\in \mathbb{Z}$, (iii) still holds after dividing out $P$ by the action of the group $\mu_{k-1}\subset S^1$ of $(k-1)$-th roots of unity. Also, the theorem has a version which applies to the case where $M$ may be non-simply connected (but then the prequantization is no longer unique). However, to avoid un-necessary complications, we restrict to the simply connected case (which, in the groupoid setting that will be adopted in the next sections, corresponds to restricting to source-simply connected groupoids).

\subsection{Outline of the paper} As we shall explain at the beginning of the next section, the prequantization problem has a very natural formulation in the world of groupoids and multiplicative forms, and, as a particular case, we also recover the notion of prequantization of symplectic groupoids introduced by Weinstein and Xu \cite{WeXu} (we allow here non-symplectic forms since, for instance, it allows us to treat also the presymplectic groupoids of \cite{CH}). Then, in this more general context, we will explain the relation between the prequantization and the integrability of a Lie algebroid naturally associated to the 2-form. The ``strategy'' is quite simple (see the end of the next section): go to the infinitesimal picture where the situation is much simpler, and then ``integrate'' back to the global picture. The supporting results for this plan are presented in 
Section \ref{integration} and Section \ref{sub-i-g}. 
More precisely, in Section \ref{integration} we recall the construction of the monodromy groupoid of an algebroid \cite{CrFe1},
and then we discuss the integrability of Lie algebroids associated to general 2-cocycles (Theorem \ref{2-coc-int}, which, we believe, is of independent interest). Then, in Section \ref{sub-i-g}, we briefly recall the reconstruction from \cite{CH} of multiplicative 2-forms out of the associated infinitesimal data, and then we prove a similar result on reconstruction of multiplicative 1-forms  (Theorem \ref{theorem-new} which, again, we believe to be of independent interest). Finally, in Section \ref{Prequantization of groupoids}, we go back to the prequantization problem, and we derive the main results on prequantization (Theorem \ref{main-theorem}, and Corollary \ref{main-theoremH} which takes care of the Hausdorff case). 

\subsection{Comparing with the existing literature} 
Our discussion here ``from prequantization to integrability'' appears implicitly in the literature in various forms. First of all, it is based on the fundamental relation between principal bundles and transitive groupoids \cite{MK}. Secondly, it is strongly related to the central extensions arising from prequantization that are described in 
Brylinski's book \cite{Br} (see also the references therein). To compare, let us point out that \cite{Br} writes the story using extensions of groups only (and not algebroids or groupoids). Of course, this forces the use of infinite dimensional Lie groups and overlooks some of the structure that is at the heart of the problem. However, like us, \cite{Br} does find a way to extend the discussion to the case where the 2-form is not integral (the ``non-integrable case''), with the remark that one has to allow certain infinite dimensional groups that are not Lie. Clearly, this is as close as \cite{Br} gets to the ``integrability problem''. 
On the other hand, ``The converse''-part of our discussion seems to be less known
(although, based on experience, I expect a large number of reactions proving the contrary).
Note that in this direction too, there are quite a few striking similarities with \cite{Br}, 
although they need (and deserve) further clarification. Refering to Remark \ref{remark-Br} for such similarities, let us point out here that 
the construction of the central extension of the loop group (section 6.4 in \cite{Br}) is completely analogous to our construction of the prequantizing bundle, but one level higher (2-homotopies instead of 1-homotopies). 
Such similarities cannot be accidental, and we believe they are just some of the small pieces of a whole (more geometric) picture based on prequantization of groupoids with a 3-form background, gerbes over such groupoids etc, picture that would be not just ``an extension to groupoids'', but also clarifying for the classical constructions of \cite{Br}, and also needed for the further geometric study. Such a picture would probably benefit also from the classification results of \cite{Mo}.

Finally, the notion of prequantization of {\it symplectic groupoids} was introduced by A. Weinstein and P. Xu in \cite{WeXu}.
In this context, our Theorem \ref{main-theorem} generalizes one of the theorems in \cite{CZ} (however, we had to find a quite different proof
which also works in the non-symplectic case, and the essential step is Theorem \ref{theorem-new}). We also mention here
the work of C. Laurent-Gengoux and P. Xu on prequantization of presymplectic groupoids \cite{LX}, with the mentionm that our work should be viewed as complementary, and based on a very different approach:
in contrast with \cite{LX}, we look at the infinitesimal picture (and this is really the main idea), and then use an integration 
step \cite{CrFe1, CH} to go back to the global picture. For this reason, the results we obtain are quite general
(in particular, they do not require any Hausdorfness assumption), and the path spaces become naturally part of the picture. 

Finally, this paper is also closely related to a paper with the same title, same author, etc, but dated 2002. If you see it, then please let me know (it might help to recover my laptop!). 
\\

{\bf Conventions:} We would like to recall here the standard convention on Lie groupoids (see e.g. \cite{CrFe1}): it is important to allow groupoids whose manifold of arrows 
is non-Hausdorff. Important examples come from foliation theory, from integrating bundles of Lie algebras, Poisson geometry, and from integrating infinitesimal Lie algebra actions on manifolds. However, the base manifold (of objects) and the $s$-fibers (consisting of arrows with fixed source) are always assumed to be Hausdorff, second countable manifolds. In this paper we also make the overall hypothesis that all manifolds (including $s$-fibers of groupoids) are connected. Also, in this paper $S^1= \mathbb{R}/\mathbb{Z}$ (i.e., the exponential map or the first Chern class will not contain the factor $2\pi \sqrt{-1}$).

\section{From classical to ``oids''}
\label{to-oids}

In this section we re-write the classical prequantization problem in terms of groupoids. This provides not only a more general framework (hence also more general theorems), but also interesting re-interpretations of some of the classical notions.

\subsection{Passing to groupoids} In what follows, for a Lie groupoid $\mathcal{G}$ over $M$, we denote by $s, t: \mathcal{G}\rmap M$ the source/target maps, and
we denote by $\mathcal{G}_2$ the space of ``composable arrows'', i.e. pairs $(g, h)$ with $s(g)= t(h)$ (see also the convention at the end of the introduction). We will denote by $u: M\rmap G$ the embedding that associates to $x$ the identity element at $x$. We say that $G$ is $s$-simply connected if the s-fibers $s^{-1}(x)$ with $x\in M$ are connected and simply connected.

Implicit in our discussion of the classical case are the following basic examples. First of all, 
the pair groupoid $G(M)= M\times M$ of a manifold $M$ is the groupoid over $M$ which has exactly one arrow between
each two points (hence source/target are the projections, and the multiplication is $(x, y)(y, z)= (x, z)$). Next, for
a principal $S^1$-bundle $P$ over $M$, the gauge groupoid of $P$, $\mathcal{G}auge(P)= G(P)/S^1$, 
is the quotient of the pair groupoid of $P$ by the action of the circle. Finally, the trivial bundle of Lie groups with fiber $S^1$,
$S^{1}_{M}= M\times S^1$, is a groupoid: the source and the target are the projection maps, and the multiplication is $(x, z_1)(x, z_2)= (x, z_2)$ (more generally, any bundle of Lie groups over $M$ 
can be viewed as a groupoid over $M$).


These three groupoids fit into a short exact sequence of Lie groupoids $S^{1}_{M}\rmap \mathcal{G}auge(P)\rmap \mathcal{G}_{M}$, and this  situation is formalized by the following definition. 

\begin{definition}
\label{ext-groupoids}
Let $G$ be a Lie groupoid over $M$, and let $\mathcal{S}$ be a bundle of abelian Lie groups.
An extension of $G$ by $\mathcal{S}$ is a short exact sequence of Lie groupoids over $M$
\[  \mathcal{S}\stackrel{i}{\rmap} \tilde{G}\stackrel{\pi}{\rmap} G ,\]
where $\pi$ is a submersion and $i$ is an embedding. When $\mathcal{S}= S^{1}_{M}$, we talk about extensions by $S^1$.
Such an extension is called central if $\gamma g= g\gamma$ for all $\gamma\in S^{1}_{M}$, $g\in \tilde{G}$
for which the products are defined.
\end{definition}

Note that, for an extension of $G$ by $S^1$, one has an induced action: for $g\in G$, and $\gamma\in S^{1}_{M}$
with $t(g)= s(\gamma)$, the expression $\tilde{g}\gamma\tilde{g}^{-1}\in S^{1}_{M}$
for $\tilde{g}\in \pi^{-1}(g)$ only depends on $g$ and $\gamma$, and this defines an action of $G$ on the bundle of groups $S^1_{M}$.
Hence, saying that an extension by $S^1$ is central is equivalent to saying that $ker(\pi)$ is trivial not only as a bundle of groups, but also as a representation of $G$. Equivalently, the induced right and left $S^1$-actions on $\tilde{G}$ coincide, and there
is no confusion when referring to $\tilde{G}$ as a principal $S^1$-bundle.

Next, intimately related to prequantization is the notion of multiplicative differential forms.
Recall that, for a groupoid $\tilde{G}$, $\tilde{G}_2$ is the space of pairs of composable arrows. There are canonical maps
\[ m, pr_1, pr_2: \tilde{G}_2\rmap \tilde{G},\]
where $m$ is the multiplication, and $pr_{i}: \tilde{G}_2\rmap \tilde{G}$ are the projections.

\begin{definition} A differential form $\omega$ on $\tilde{G}$ is called multiplicative if $m^*\omega= pr_{1}^{*}\omega+ pr_{2}^{*}\omega$. 
\end{definition}

Back to the motivating example, we have:

\begin{lemma} Modulo isomorphisms, there is a 1-1 correspondence between principal $S^1$-bundle over $M$ and central extensions
of $G(M)$ by $S^1$, which associates to a bundle $P$ the extension 
\[ S^{1}_{M}\rmap \mathcal{G}auge(P) \stackrel{\pi}{\rmap} G(M) .\]
Moreover, given the bundle $\pi_{P}: P\rmap M$, there are 1-1 correspondences $\omega\mapsto \tilde{\omega}$ and $\theta\mapsto \bar{\theta}$ between 
\begin{eqnarray*}
\text{(closed)\ 2-forms\ $\omega$\ on\ $M$} & \longleftrightarrow & \text{(closed)\ multiplicative\ 2-forms\ $\tilde{\omega}$\ on\ $G(M)$}\\
\text{$S^1$-invariant\ 1-forms\ $\theta$\ on\ $P$} & \longleftrightarrow & \text{multiplicative\ 1-forms\ $\bar{\theta}$\ on\ $\mathcal{G}auge(P)$}
\end{eqnarray*}
and the formula $\pi_{P}^{*}\omega= d\theta$ is equivalent to $\pi^{*}\tilde{\omega}= d\bar{\theta}$.
\end{lemma}


The correspondences above are easy to describe: for $\omega\in \Omega^2(M)$,
$\tilde{\omega}:= pr_{1}^{*}\omega- pr_{2}^{*}\omega$,  
while for $\theta\in \Omega^{1}(P)^{S^1}$, $\tilde{\theta}= pr_{1}^{*}\theta- pr_{2}^{*}\theta$ descends to a multiplicative 1-form on the gauge groupoid $\mathcal{G}auge(P)$,
and $\bar{\theta}$ is the resulting form. The remaining details for the proof of the lemma are an instructive exercise.
With these in mind, we have the following definition (see also \cite{WeXu, LX}).

\begin{definition}\label{preq-g} Let $G$ be a Lie groupoid over $M$, and let $\omega$ be a closed multiplicative 2-form on $G$. A prequantization of $(G, \omega)$ is a central extensions of Lie groupoids:
\[  S^{1}_{M}\rmap \tilde{G}\stackrel{\pi}{\rmap} G ,\]
and a multiplicative 1-form $\theta\in \Omega^1(\mathcal{G})$ which is a connection form for the principal $S^1$-bundle $\pi: \tilde{G}\rmap G$ and which satisfies $d\theta= \pi^*\omega$.
\end{definition}

\subsection{The infinitesimal analogue}
\label{The infinitesimal analogue} We now describe the infinitesimal picture underlying the situation that appears in the previous definition. Hence, let $A$ be a Lie algebroid over $M$ (the infinitesimal counterpart of the groupoid $G$). 

\begin{definition}
An extension of $A$ by $\mathbb{R}$ is a short exact sequence of algebroids over $M$
\[ \mathbb{R}_{M}\stackrel{i}{\rmap} \tilde{A}\stackrel{\pi}{\rmap} A ,\]
where $\mathbb{R}_{M}$ is the trivial real line bundle over $M$, viewed as a bundle of abelian Lie algebras (hence also 
as an algebroid with zero anchor).  
Such an extension is called central if
$[\tilde{\alpha}, i(1)]= 0$ for all $\tilde{\alpha}\in \tilde{A}$.
\end{definition}

Central extensions are well-known in various other contexts (e.g. groups, Lie algebras etc), as is the fact that they are related to (and even classified by) 2-cocycles. To review this in the case of Lie algebroids, recall that the ``algebroid DeRham complex'' associated to
$A$, $(C(A), d_{A})$, is defined as follows: $C^p(A)= \Gamma(\Lambda^pA^*)$ consists of
antisymmetric $C^{\infty}(M)$-multilinear maps
\[ c: \underbrace{\Gamma(A)\otimes \ldots \otimes \Gamma(A)}_{p\ \text{times}}\rmap C^{\infty}(M),\]
with the differential
\begin{eqnarray*}
d_A(c)(\alpha_1, \ldots , \alpha_{p+1}) & = & \sum_{i<j}
(-1)^{i+j}c([\alpha_i, \alpha_j], \alpha_1, \ldots ,
\hat{\alpha_i}, \ldots ,
\hat{\alpha_j}, \ldots \alpha_{p+1}) \\
 & + & \sum_{i=1}^{p+1}(-1)^{i+1}
\mathcal{L}_{\rho(\alpha_i)}(c(\alpha_1, \ldots, \hat{\alpha_i},
\ldots , \alpha_{p+1})) .
\end{eqnarray*}
A 2-cocycle on $A$ will be an element $c\in C^2(A)$ such that
$d_{A}(c)= 0$. Any such cocycle $c$ induces a central extension of $A$ by $\mathbb{R}$,
$\mathbb{R}_{M}\stackrel{i}{\rmap} A_{c}\stackrel{\pi}{\rmap} A$, 
where 
\[ A_c= A\oplus \mathbb{R},\] 
with the bracket:
\[ [(\alpha, f), (\beta, g)]= ([\alpha, \beta], \mathcal{L}_{\rho(\alpha)}(g)- \mathcal{L}_{\rho(\beta)}(f)+ c(\alpha, \beta)),\]
and with the anchor $(\alpha, f)\mapsto \rho(\alpha)$. Actually, one can check that $A_c$ is an algebroid if and only if $c$ is a cocycle.

What makes the infinitesimal picture useful is that, viewing a 2-cocycle $c\in C^2(A)$ as a replacement of multiplicative 2-forms on groupoids, the
infinitesimal analogue of the prequantization data (i.e. of the Definition \ref{preq-g}) comes for free: the projection on the
second component,
\[ l_{c}(\alpha, \lambda):= \lambda,\]
defines an element $l_{c}\in C^1(A)$ with the properties:
\[ l_c(i(1))= 1, \ d_{A}(l_c)= \pi^*(c). \]

\subsection{From groupoids to algebroids}
\label{From groupoids to algebroids}
We now explain the passing from central extensions and multiplicative 2-forms on Lie groupoids to central extensions and 2-cocycles of Lie algebroids. First of all, any Lie groupoid $G$ has an associated Lie algebroid $A= A(G)$, constructed in the same way that one constructs the Lie algebra of a Lie group: as a vector bundle, $A= Ker(ds)|_{M}$, where the restriction is with respect to the embedding $u: M\rmap G$ which associates to $x\in M$ the identity arrow at $x$ (hence $A_x$ is the tangent space at $u(x)$ of $s^{-1}(x)$). Next, $\Gamma(A)$ is isomorphic to the space $\mathcal{X}^{s}(M)^{inv}$ of vector fields tangent to the $s$-fibers which are invariant under right translations. This induces a bracket on $\Gamma(A)$, while $(dt)_{u(x)}: T_{u(x)}s^{-1}(x)\rmap T_xM$ defines the anchor of $A$. Given a Lie algebroid $A$, we say that $A$ is integrable if it is isomorphic to the Lie algebroid of a Lie groupoids. In general, Lie algebroids may fail to be integrable. However, if $A$ is integrable then, similar to Lie algebras, there exists a unique $s$-simply connected (i.e. with connected and simply connected $s$-fibers) Lie groupoid integrating $A$.

Next, any closed multiplicative 2-form $\omega$ on $G$ induces (by restriction) a 2-cocycle $c_{\omega}$ on $A= A(G)$:
\[ c_{\omega}(\alpha, \beta)= \omega(\alpha, \beta),\ \ \ \text{for}\ \alpha, \beta\in A_{x}\subset T_{u(x)}G .\]
Similarly, starting with a central extension of $G$, and passing to Lie algebroids, one obtains a central extension of $A$.

\begin{lemma}
\label{first-step}
Given a closed multiplicative 2-form on $\omega$ on $G$, then any prequantizing groupoid $\tilde{G}$ of $(G, \omega)$ integrates $A_{c_{\omega}}$.
\end{lemma}

\begin{proof} Assume that $\tilde{G}$ is a prequantization of $(G, \omega)$ with connection 1-form $\theta$. Passing to algebroids, we obtain an extension $\mathbb{R}\rmap \tilde{A}\stackrel{\pi}{\rmap} A$,
where $\tilde{A}$ is the Lie algebroid of $\tilde{G}$. Moreover, the 1-form $\theta$ defines, by restriction to $A$, an element $l\in C^1(A)$, and the
prequantization conditions on $\theta$ translate into $d_{\tilde{A}}(l)= \pi^{*}(c_{\omega})$ and $l(1)= 1$.
But this precisely means that the map
\[ \tilde{A}\rmap A_{c_{\omega}},\ \alpha\mapsto (\pi(\alpha), l(\alpha)) \]
defines an isomorphism of extensions, sending $l$ to the canonical $l_{c_{\omega}}\in C^1(A_{c_{\omega}})$. In particular, $\tilde{A}$ is isomorphic to $A_{c_{\omega}}$, hence $\tilde{G}$ integrates $A_{c_{\omega}}$.
\end{proof}

\subsection{Prequantizing: the strategy} With Lemma \ref{first-step} in mind, we have the following strategy for reconstructing the prequantization
of $(G, \omega)$: we consider the algebroid 2-cocycle $c$ which is the restriction of $\omega$ to $A$, we integrate $A_c$ to a groupoid $G_{\omega}$,
we recover the connection 1-form from its infinitesimal data, and then we eventually ``correct'' $G_{\omega}$ so that it becomes an extension by $\mathbb{S}^1$ (the last step is needed: for instance, $G_{\omega}$ may be an extensions by $\mathbb{R}$- and then, of course, we will divide out by the action of $\mathbb{Z}$). These are three steps, hence three sections still to come.

\section{Integration}
\label{integration}

\subsection{Integrating Lie algebroids}
\label{sub-i-g2}
We now briefly recall the construction of the s-simply connected groupoid associated to an algebroid \cite{CrFe1}.
So, let $A$ be an algebroid over $M$. We consider the Banach manifold
$\tilde{P}(A)$ consisting of paths $a:I\rmap A$ from the unit interval $I$ into $A$,
of class $C^1$, whose base path $\gamma= p\circ a:I\rmap M$ is of class $C^2$ (where  $p:A \rmap M$ is the natural projection).
Denote by $P(A)$ the submanifold consisting of $A$-paths, i.e. those $a\in \tilde{P}(A)$ which satisfy the equation $\rho\circ a= \frac{d}{dt}\gamma$.
Next, we need the notion of ``$A$-homotopies with fixed end-points''. To describe this, we choose a connection on the vector bundle $A$, and we consider its $A$-torsion $T_{\nabla}(\alpha, \beta)= \nabla_{\rho(\alpha)}(\beta)-  \nabla_{\rho(\beta)}(\alpha)- [\alpha, \beta]$. Using the connection, for any path $a= a(t): I\rmap A$ we can talk about its derivative $\partial_{t}a:= \nabla_{\frac{d\gamma}{dt}}(a): I\rmap A$, where $\gamma$ is the base path of $a$. Given a family $a_{\epsilon}$ ($\epsilon\in I$) of $A$-paths with base paths $\gamma_{\epsilon}= \gamma_{\epsilon}(t)$, we consider $b= b(\epsilon, t)$ (sitting over $\gamma_{\epsilon}(t)$) which are defined as the solution of the equation
\begin{equation}
\label{A-homot}
\partial_tb- \partial_{\epsilon}a= T_{\nabla}(a,
b),\ \ b(\epsilon, 0)= 0 .
\end{equation}
We will
say that the family of $A$-paths $a_{\epsilon}$ is an $A$-homotopy
(with fixed end points) if $b(\epsilon,
1)= 0$ (note that $b(\epsilon, 1)$ only depends on $a_{\epsilon}$
and not on the choice of the connection!). The previous equation should be viewed as an analogue of the equation
``$\frac{d}{d\epsilon}\frac{d}{dt}= \frac{d}{dt}\frac{d}{d\epsilon}$''.
Thinking of $A$ as a generalized tangent bundle, an $A$-path $a$ can be thought of as a path $\gamma$ in $M$ together with an ``$A$-derivative''
(the path $a$ sitting above $\gamma$) which gives back the true derivative of $\gamma$ after applying the anchor. The role of the previous equation is: starting from $a$ (the $A$-derivative of $\gamma$ with respect to $t$) find the $A$-derivative of $\gamma$ with respect to $\epsilon$ (i.e. $b$). This can be made more precise if $A$ comes from an $s$-simply connected Lie groupoid $G$: then any
path $g$ in $G$ which stays in an $s$-fiber $s^{-1}(x)$ of $G$ induces an $A$-path
$a= D^{R}g$ by taking derivatives and translating back to an identity element using right translation:
\begin{equation}
\label{bij-paths}
a(t)= (dR_{g(t)^{-1}})_{g(t)} \dot{g}(t) \ .
\end{equation}
This defines a bijection between $A$-paths and paths in $G$ which stay in a single $s$-fiber of $G$ and start at an identity element.
Similarly, families $a_{\epsilon}$ are identified with families $g_{\epsilon}$ in $G$ which stay in a single $s$-fiber, and the solution $b$ comes from right translation of the derivative $\frac{dg}{d\epsilon}$: $b= D_{\epsilon}^{R}(g)$.

Denote by $\sim$ the $A$-homotopy equivalence relation on $P(A)$, and define
the monodromy groupoid of $A$,
\[ G(A):= P(A)/\sim .\]
It is a groupoid with the source (resp. target) map obtained by
taking the starting (resp. ending) point of the base paths, and
the multiplication is defined by concatenation of paths \cite{CrFe1}. Note that
if $A$ is a Lie algebra (i.e. $M$ is a point), then $G(A)$ is the
(unique) simply connected Lie group integrating $A$, as
constructed in \cite{DuKo}.

Finally, we recall the notion of monodromy. Let $x\in M$. We
denote by $L\subset M$ the orbit of $A$ through $x$. It consists
of those points $y$ with the property that there exists an
$A$-path $a$ whose base bath joins $x$ and $y$, and it is also the
maximal integral submanifold integrating the singular distribution
$\rho(A)\subset TM$. We also denote by $\mathfrak{g}_{x}(A)=
Ker(\rho_x)$ the isotropy Lie algebra of $A$ at $x$. The bracket of $A$
restricts to a Lie algebra bracket on $\mathfrak{g}_{x}(A)$:
if $\alpha(x), \beta(x)\in \mathfrak{g}_x(A)$, the expression
$[\alpha, \beta](x)$ only depends on $\alpha(x)$ and $\beta(x)$. 
With the previous notations, $G(\mathfrak{g}_x(A))$ is the
simply connected Lie group integrating $\mathfrak{g}_x(A)$.
The monodromy map at $x$,
\[ \partial_x: \pi_2(L, x)\rmap G(\mathfrak{g}_x(A)) \]
is defined as follows. Given $\gamma= \gamma(\epsilon, t): I^2\rmap L$ with $\gamma(\partial I^2)= x$, one can always find an $A$-homotopy $a_{\epsilon}= a_{\epsilon}(t)$ of $A$-paths above $\gamma$, with $a_0= 0$. Then $a_1= a_1(t)$ will be a path in $\mathfrak{g}_x(A)$, and 
\[ \partial_x([\gamma]):= [a_1] .\]
The monodromy group at $x$, denoted $\mathcal{N}_x(A)$, is defined
as the image of $\partial_x$. These groups arise naturally when comparing
the isotropy group $G(A)_x= s^{-1}(x)\cap t^{-1}(x)$ of $G(A)$ with $G(\mathfrak{g}_x(A))$. 
The two groups are made out of the same paths, but using (slightly) different equivalence relations. It is not difficult to see \cite{CrFe1} that the connected component of $G(A)_x$ containing the identity is
\begin{equation}
\label{id-comp}
G(A)_{x}^{0}=
G(\mathfrak{g}_x(A))/\mathcal{N}_x(A).
\end{equation}

As for $A$-paths and homotopies, $\partial_x$ and $\mathcal{N}_x(A)$, too,
become more transparent if $A$ comes from an $s$-simply connected Lie groupoid $G$.
First of all, $G(x, x)=
s^{-1}(x)\cap t^{-1}(x)$ will be a Lie group that integrates
$\mathfrak{g}_x(A)$, hence $\pi_1(G(x, x))$ sits inside the simply
connected Lie group $G(\mathfrak{g}_x(A))$ integrating
$\mathfrak{g}_x(A)$, and this is precisely $\mathcal{N}_x(A)$. On
the other hand, $t: s^{-1}(x)\rmap L$ is a principal $G(x,
x)$-bundle, and then $\partial$ will be precisely the boundary map
of the induced homotopy long exact sequence:
\begin{equation}
\label{ses-}
\ldots \rmap \pi_2(s^{-1}(x))\rmap
\pi_2(L)\stackrel{\partial}{\rmap} \pi_1(G(x, x))(\subset
G(\mathfrak{g}_x(A))).
\end{equation}

We can now recall the main result of \cite{CrFe1}:

\begin{theorem}\label{th-CF} For any Lie algebroid $A$, the orbits of the equivalence relation $\sim$ ($A$-homotopy) define a regular 
foliation on $P(A)$ of finite codimension equal to $dim(M)+ rank(A)$, $G(A)$ is a topological groupoid, and the following are equivalent:
\begin{enumerate}
\item[(i)] $A$ is integrable.
\item[(ii)] $G(A)$ is smooth.
\item[(iii)] The monodromy groups $\mathcal{N}_{x}(A)$ are locally uniformly (with respect to $x\in M$) discrete.
\end{enumerate}
Moreover, in this case $G(A)$ will be the unique $s$-simply connected Lie groupoid integrating $A$.
\end{theorem}

For the study of the smoothness of $G(A)$ (including the proof of the theorem) and of the geometric structures on it, 
it is important to realize that
the $A$-homotopy equivalence classes on $P(A)$ can be realized as the orbits of a Lie algebra action. 
The Lie algebra, denoted $\mathfrak{g}$, consists of time dependent sections $\eta_{t}$ ($t\in [0, 1]$) of $A$, with $\eta_0= \eta_1= 0$, and the infinitesimal action on $P(A)$ translates into a Lie algebra map
\begin{equation}
\label{vctrfld}
\mathfrak{g}\ni \eta\mapsto X_{\eta}\in \mathcal{X}(P(A)) .
\end{equation}
To define $X_{\eta}$, we fix $a_0\in P(A)$, and we construct  the flow
$a_{\epsilon}= \Phi_{X_{\eta}}^{\epsilon}(a_0)$ in such a way that $a_{\epsilon}$
are paths above $\gamma_{\epsilon}(t)= \Phi_{\rho(\eta_t)}^{\epsilon}\gamma_{0}(t)$, where $\gamma_{0}$ is the
base path of $a_0$, and $\Phi_{\rho(\eta_t)}^{\epsilon}$ is the flow of the vector field $\rho(\eta_t)$.
We choose a time dependent section $\xi_{0}$ of $A$ with $\xi_{0}(t, \gamma_0(t))= a_0(t)$,
and we consider the $(\epsilon, t)$-dependent section of $A$,
$\xi= \xi(\epsilon, t)$, solution of
\begin{equation}
\label{def-xi}
\frac{d\xi}{d\epsilon}- \frac{d\eta}{dt}= [\xi, \eta], \; \; \xi(0, t)= \xi_0(t).
\end{equation}
Then $a_{\epsilon}(t)= \xi_{\epsilon}(t, \gamma_{\epsilon}(t))$.
This defines the desired vector fields $X_{\eta}$, hence the action of $\mathfrak{g}$.

\subsection{The case of 2-cocycles}
\label{The case of 2-cocycles}
Assume now that $A$ is the Lie algebroid of an $s$-simply
connected Lie groupoid $G$, $c\in C^2(A)$ is a 2-cocycle, and we
look at the integrability of the associated algebroid $A_c$ (see subsection \ref{The infinitesimal analogue}). For
$x\in M$, we use right translations to extend $c\in
\Gamma(\Lambda^2A^*)$ to a 2-form on the
$s$-fiber $s^{-1}(x)$, call it $\omega_{c}^{x}$. In other words, for $g\in s^{-1}(x)$ and 
$X_g, Y_g\in T_gs^{-1}(x)$, one uses the right translation $R_{g^{-1}}: s^{-1}(x)\rmap s^{-1}(y)$
($y$ is the target of $g$) and its differential at $g$, and
\begin{equation}
\label{omega-c-fr} 
\omega_{c}^{x}(X_{g}, Y_{g})= c_y((dR)_{g^{-1}}(X_g), (dR)_{g^{-1}}(Y_g)).
\end{equation}
We define the period group and the structural group of $c$ at $x$ as:
\[ \mathcal{P}_x(c):= Per(\omega_{c}^{x})\subset \mathbb{R}, \ \mathcal{S}_{x}(c):= \mathbb{R}/\mathcal{P}_x(c).\]
Varying $x\in M$, we obtain two (possibly non-smooth) bundles of
groups over $M$, denoted by $\mathcal{P}(c)$ and $\mathcal{S}(c)$,
respectively. There is no ambiguity when talking about the smoothness of these bundles: we ask $\mathcal{P}(c)$ to be a smooth sub-bundle of $\mathbb{R}_{M}$, and we ask the projection from $\mathbb{R}_{M}$ into $\mathcal{S}(c)$ to be a submersion.

\begin{theorem}\label{2-coc-int} If $A$ is the Lie algebroid associated to the $s$-simply
connected Lie groupoid $G$, $c\in C^2(A)$ is a 2-cocycle, and
$A_c$ is the associated algebroid, then there is an extension of topological groupoids
\begin{equation}
\mathcal{S}(c)\rmap G(A_c)\stackrel{\pi}{\rmap} G
\end{equation}
and the following are equivalent:
\begin{enumerate}
\item[(i)] $A_c$ is integrable.
\item[(ii)] $G(A_c)$ is smooth.
\item[(iii)] $\mathcal{P}(c)$ is smooth.
\item[(iv)] $\mathcal{S}(c)$ is smooth.
\end{enumerate}
\end{theorem}

\begin{proof} 
We may assume that $G= G(A)$ and let $\pi$ be the groupoid map induced by
the projection $A_c\rmap A$. We want to show that $Ker(\pi_x)$ can
be identified with $\mathcal{S}_x(c)$. Let $x\in M$ and let $L$ be the
leaf through $x$. We denote by $\tilde{\partial}: \pi_2(L)\rmap G(\mathfrak{g}_x(A_c))$ the monodromy map of
$A_c$ and by $\partial$ the one of $A$. First of all, we have
\begin{equation}
\label{product}
G(\mathfrak{g}_x(A_c))= G(\mathfrak{g}_x(A))\times
\mathbb{R}
\end{equation} 
which follows from the similar relation at the Lie algebra level. We claim that, with respect to this decomposition,
$\partial$ is the first component of $\tilde{\partial}$, i.e. 
\begin{equation}
\label{to-prove-1}
\tilde{\partial}_x(\gamma)=
(\partial_x(\gamma), r(\gamma))
\end{equation}
for some $r: \pi_2(L)\rmap \mathbb{R}$ that we are going to make more explicit. 
Due to the relation
between the isotropy group of $G(A)$ and monodromy (see equation
(\ref{id-comp})), we deduce that $Ker(\pi_x)$ is
$\mathbb{R}/\{ r(\gamma): \partial(\gamma)= 1\}$. 
On the other hand, using the description of $\partial$ as the
boundary map of the homotopy long exact sequence (see
(\ref{ses-})), this will be the quotient of $\mathbb{R}$ by 
$\{ r(t\circ g): [g]\in \pi_2(s^{-1}(x))\}$.
Our second claim
is that, for $g: I^2\rmap s^{-1}(x)$ representing a second
homotopy class in $s^{-1}(x)$,
\begin{equation}
\label{to-prove-2} 
r(t\circ g)= \int_{g}\omega_{c}^{x} .
\end{equation}
Clearly, these imply that $Ker(\pi_x)= \mathcal{S}_x(c)$.

To prove the two claims we made (equations (\ref{to-prove-1}) and (\ref{to-prove-2})), we first rewrite the homotopy of $A_c$-paths in terms of
the homotopy of $A$-paths. It is clear that 
$A_c$-paths are pairs $\tilde{a}= (a, f)$ where $a$ is an $A$-path and $f: I\rmap \mathbb{R}$. 
Next, the condition that a family $\tilde{a}_{\epsilon}= (a_{\epsilon}, f_{\epsilon})$ is an $A_{c}$-homotopy
brakes into two conditions
(depending on the components we look at): first of all, $a_{\epsilon}$ must be an $A$-homotopy, i.e.
the solution $b$ of (\ref{A-homot}) satisfies $b(\epsilon, 1)= 0$. Secondly,
the solution $h:I\rmap \mathbb{R}$ of the equation
\begin{equation}
\label{hom-scalar} 
\frac{df}{d\epsilon}- \frac{dh}{dt}= c(a, b) ,\ h(\epsilon, 0)= 0
\end{equation}
must satisfy $h(\epsilon, 1)= 0$. Computing $h(\epsilon, 1)$ by integrating the previous equation, we conclude that 
$\tilde{a}_{\epsilon}= (a_{\epsilon}, f_{\epsilon})$ is an $A_{c}$-homotopy if and only if $a_{\epsilon}$ is an $A$-homotopy
and 
\begin{equation}
\label{to-integrate} 
\frac{d}{d\epsilon} \int_{0}^{1} f(\epsilon, t) dt= \int_{0}^{1} c(a, b) dt .
\end{equation}
This discussion immediately implies that $\partial$ is the first component of $\tilde{\partial}$.
Indeed, choosing $\tilde{a}_{\epsilon}= (a_{\epsilon}, f_{\epsilon})$ as in the definition of $\tilde{\partial}(\gamma)$ (i.e. an $A_c$-homotopy above $\gamma$, with $\tilde{a}_{0}= 0$),   
$a_{\epsilon}$ will satisfy the similar 
conditions appearing in the definition of $\partial(\gamma)$. Hence, with respect to the
decomposition (\ref{product}), 
 $\tilde{\partial}(\gamma)= ([a_1], [f_1])= (\partial(\gamma), [f_1])$. Moreover, since $G(\mathbb{R})= \mathbb{R}$ identifies the $\mathbb{R}$-homotopy class of an $\mathbb{R}$-path $f_1= f_{1}(\epsilon)$ with is integral $\int_{0}^{1}f_1d\epsilon \in \mathbb{R}$, and the last integral can be computed by integrating (\ref{to-integrate}), we deduce that the second component of $\tilde{\partial}(\gamma)$ 
(see \ref{to-prove-1}) is
\begin{equation}
\label{last-fr} 
r(\gamma)= \int_{0}^{1}\int_{0}^{1} c(a, b) dt d\epsilon .
\end{equation}
In this formula, we have the freedom of choosing any $a$ and $b$ that satisfy the conditions above: they sit over $\gamma$, 
and they are as in the definition of $\partial(\gamma)$. We can now prove equation (\ref{to-prove-2}). Starting with $g= g_{\epsilon}(t)$, we choose as $A$-homotopy the family $a_{\epsilon}= D^{R}_{t}(g)$ induced by $g_{\epsilon}$ by taking derivatives with respect to $t$ and using right translations (see equation (\ref{bij-paths})), and then the $b$'s will be obtained similarly by interchanging the role of $\epsilon$ and $t$. On the other hand, exactly the same right translations appear in the definition of  
$\omega_{c}^{x}(\frac{dg}{dt}, \frac{dg}{d\epsilon})$ (see equation (\ref{omega-c-fr})), and this translates into 
\[ \omega_{c}^{x}(\frac{dg}{dt}, \frac{dg}{d\epsilon})= c(a, b) .\]
Using this in (\ref{last-fr}) $r(\gamma)$, we obtain (\ref{to-prove-2}).

We are left with proving the equivalence of (i)-(iv). The
implications (ii)$\Longrightarrow$(iii)$\Longrightarrow$(iv) are
easy (all the spaces to be shown to be smooth are inverse images of
submanifolds by submersions), while (i)$\Longrightarrow$(ii) is
part of Theorem \ref{th-CF} applied to $A_c$. To prove
(iv)$\Longrightarrow$(i) we use again Theorem \ref{th-CF}
applied to $A_c$, and we check that $Im(\tilde{\partial}_x)$ is
locally uniformly discrete. Hence, assume that
$\tilde{\partial}_x(\gamma_n)$ converge to $1$, and we want to
prove that all these elements equal $1$ for $n$ large enough. From
(\ref{to-prove-1}) we deduce that $\partial(\gamma_n)$ converges
to $1$, and similarly $r(\gamma_n)$. Since $A$ is integrable,
$\partial(\gamma_n)= 1$ for $n$ large enough, and then, as proven
above, $r(\gamma_n)\in Per(\omega_{c}^{x})= \mathcal{P}_{x}(c)$.
But the projection $\mathcal{P}\rmap M$ is a local diffeomorphism
because its fibers are at most countable subgroups of $\mathbb{R}$
(second homotopy groups of manifold have this property), hence $r(\gamma_n)= 1$ for $n$ large enough.
\end{proof}

\begin{remark}\rm \ Note that the general arguments on $A$-homotopies given in the previous proof give the following simplified description
of $G(A_c)$: any $A_c$-path is homotopic to one of type $(a, r)$ where $a$ is an $A$-path and $r\in \mathbb{R}$. Indeed, each $(a, f)$ is homotopic to $(a, \int_{0}^{1} f)$ by the homotopy with $a_{\epsilon}= a$ (the corresponding $b$ will be $0$),
and $f_{\epsilon}= \epsilon f+ (1-\epsilon)\int_{0}^{1} f$.

Also, two such paths
$(a_0, r_0)$ and $(a_1, r_1)$ with $r_i\in \mathbb{R}$ are homotopic if and only if there is a homotopy $a$ between $a_0$ and $a_1$ such that
$r_1- r_0= \int_{I^2} c(a, b) d\epsilon dt$. 
In one direction, if $(a_{\epsilon}, f_{\epsilon})$ is a homotopy between these two elements, one just integrates the equation (\ref{to-integrate}). Conversely, if $a_{\epsilon}$ is an $A$-homotopy which satisfies the last equation, then 
$(a_{\epsilon}, \epsilon r_0+ (1-\epsilon) r_1)$ will satisfy (\ref{to-integrate}), i.e. defines an $A_c$-homotopy.

In conclusion, $G(A_c)$ is the quotient of $P(A)\times \mathbb{R}$ by the equivalence relation described by the formula above. 
In the case where $A= TM$ (hence $c= \omega$ is a closed 2-form on $M$) we obtain a transitive groupoid and, restricting to its s-fiber above $x_0\in M$, we obtain the description of $P_{\omega}$ in the introduction.  
\end{remark}


\section{Reconstructing multiplicative forms}
\label{sub-i-g}

\subsection{Reconstructing closed multiplicative 2-forms}
We have remarked in Section \ref{to-oids} that multiplicative 2-forms on a Lie groupoid $G$ induce 2-cocycles on the associated
Lie algebroid $A$. However, such cocycles are not the ``hole infinitesimal counterpart'' of multiplicative forms.
Finding the right answer is a bit more subtle, and it has been done in \cite{CH}. The key remark is that
any multiplicative s-form $\omega\in\Omega^2(G)$ induces a bundle map
\[ \rho^{*}_{\omega}: A\rmap T^*M,\]
by the equation
\[ \SP{\rho_{\omega}^{*}(\alpha), X}= \omega(\alpha, X) ,\]
where $\SP{\cdot, \cdot}$ denotes the pairing between cotangent and tangent vectors, and $TM$ is viewed as a subspace of $TG$
by using the embedding $u: M\rmap G$. If $\omega$ is closed, then this map satisfies the basic relations
\begin{eqnarray}
& \SP{\rho^{*}_{\omega}(\alpha),\rho(\beta)}   & =
                        -\SP{\rho^{*}_{\omega}(\beta),\rho(\alpha)}\label{mult-c-1}\\
& \rho^{*}_{\omega}([\alpha, \beta])   & =
\mathcal{L}_{\alpha}(\rho^{*}_{\omega}(\beta)) - \mathcal{L}_{\beta}(\rho^{*}_{\omega}(\alpha))+
d\SP{\rho^*_{\omega}(\alpha),\rho(\beta)} ,\label{mult-c-2}
\end{eqnarray}
for all $\alpha, \beta\in \Gamma(A)$. Note that the 2-cocycle $c_{\omega}$ associated to $\omega$ can be recovered from
$\rho^{*}_{\omega}$ by
\[ c_{\omega}(\alpha, \beta)= \SP{\rho_{\omega}^{*}(\alpha), \rho(\beta)},\]
and the previous formulas for $\rho_{\omega}^{*}$ imply that
$c_{\omega}$ is, indeed, a 2-cocycle. We have \cite{CH}:

\begin{theorem}
\label{th-CH} If $G$ is an $s$-simply connected Lie groupoid, then
the correspondence $\omega\mapsto \rho^{*}_{\omega}$ defines a
bijection
 between closed multiplicative 2-forms on $G$ and bundle maps $\rho^{*}: A\rmap T^*M$ satisfying (\ref{mult-c-1}) and (\ref{mult-c-2}).
\end{theorem}

The proof is based on the explicit construction of $G(A)$ (since $G$ is $s$-simply connected, we may assume that $G= G(A)$), and we briefly
recall the reconstruction of $\omega$ out of $\rho^*$. First of all, we consider the Liouville contact 1-form $\sigma^c$ on $T^*M$:
\[
\sigma^{c}_{\xi_x}(X_{\xi_x})=
\SP{\xi_x, (dp)_{\xi_{x}}(X_{\xi_x})} ,\ \ (\xi_x\in T_{x}^{*}M, X_x\in T_{\xi_x}T^*M)
\]
where $p: T^*M\rmap M$ is the projection. Using $\rho^*: A\rmap T^*M$, we pull
$\sigma^{c}$ back to $A$, and we denote by $\tilde{\sigma}$ the resulting form on $P(A)$. Hence
\[
\tilde{\sigma}_{a}(X_a)= \int_{0}^{1} \SP{ \rho^*(a(t)), (dp)_{a(t)}(X_a(t))} dt .
\]
It follows that, for vector fields (\ref{vctrfld}) coming from the action of the Lie algebra $\mathfrak{g}$ described in subsection \ref{sub-i-g2}, 
(i.e. the $X_{\eta}$'s described there), we have
\begin{eqnarray}
i_{X_{\eta, a}}(\tilde{\sigma}) = -\int_{0}^{1} \SP{{\rho^{*}\eta(t, \gamma(t))},\rho(a(t))} dt,
\;\; \mbox{ and }\label{e1-form}\\
\mathcal{L}_{X_{\eta}}(\tilde{\sigma})=  -
d(\int_{0}^{1} \SP{ {\rho^{*}\eta(t, \gamma(t))}, \frac{d\gamma}{dt} } dt),\label{e2-form}
\end{eqnarray}
where the right hand side comes from the differential of the function on $P(A)$ which sends an $A$-path $a$ into
$\int_{0}^{1} \SP{\rho^{*}\eta(t, \gamma(t)),\frac{d\gamma}{dt}} dt$
(see Lemma 5.2 in \cite{CH}). One deduces that $-d\tilde{\sigma}$ is a closed
two-form on $P(A)$ which is basic with respect to the action of
$\mathfrak{g}$, hence it will descend to a two from $\omega$ on
$G(A)$, and this will be the desired multiplicative form. For more
details, see \cite{CH}.

\subsection{Reconstructing multiplicative 1-forms} Next, we would like to know when a multiplicative 2-form is ``multiplicatively exact'',
i.e. is of type $d\theta$ for some multiplicative 1-form $\theta$. We will prove the following.

\begin{theorem}
\label{theorem-new} Let $G$ be an $s$-simply connected Lie
groupoid, let $\omega\in \Omega^2(G)$ be a closed multiplicative
2-form, and let $c_{\omega}\in C^2(A)$ be the induced algebroid
2-cocycle (the restriction of $\omega$ to $A$). Then $\omega$ is
multiplicatively exact if and only if $c_{\omega}$ is exact as an
algebroid cocycle. More precisely, there is a 1-1 correspondence
between
\begin{itemize}
\item $\theta\in \Omega^1(G)$ multiplicative such that $d\theta= \omega$.
\item $l\in C^{1}(A)$ such that $d_{A}(l)= c_{\omega}$.
\end{itemize}
\end{theorem}

\begin{proof} In one direction, given $\theta$, we choose $l$ to be the restriction of $\theta$ to $A$, and then, writing
$d\theta= \omega$ on elements $\alpha, \beta\in A$, after a short computation, we obtain $d_{A}(l)= c_{\omega}$. This is similar (but simpler) to Proposition 3.5. (ii) in \cite{CH}. For the converse,
we may assume, as before, that $G= G(A)$. Given $l\in C^1(A)$, we consider the function $f_{l}$ on $P(A)$ defined by
\[ f_{l}(a)= \int_{a}l= \int_{0}^{1} \SP{l(\gamma(t)), a(t)} dt .\]
If $l$ was a 1-cocycle (i.e. $d_{A}(l)= 0$), then this formula would only depend on the homotopy class of $a$ (or, equivalently, $f_{l}$ is invariant with respect to the action of $\mathfrak{g}$ on $P(A)$). What is left of this property when $l$ is no longer closed is the following
formula computing the Lie derivatives of $f_{l}$ with respect to vector fields
$X_{\eta}$ coming from $\mathfrak{g}$:
\begin{equation}
\label{nice} L_{X_{\eta}}(f_l)(a)= -\int_{0}^{1}
(d_{A}l)(a(t), \eta(t, \gamma(t)) dt .
\end{equation}
To prove this, we make use of the description of the flow of $X_{\eta}$ (subsection \ref{sub-i-g2}), hence we choose $\xi= \xi(\epsilon, t)$ as in equation (\ref{def-xi}). We see that the left hand side of (\ref{nice})
is the integral, with respect to $t$, of
\[
\begin{split}
\left.\frac{d}{d\epsilon}\right |_{\epsilon= 0} \SP{l(\gamma_{\epsilon}(t)), \xi_{\epsilon}(t, \gamma_{\epsilon}(t))}= \\
= \SP{l, \left.\frac{d}{d\epsilon}\right |_{\epsilon= 0} \xi_{\epsilon}}+ L_{\frac{d\gamma}{d\epsilon}}(\SP{l, \xi}) \ \ \text{(at\ the\ point}\ x= \gamma(t))\\
= \SP{l, \frac{d\eta}{dt}}+ \SP{l, [\xi, \eta]}+ \mathcal{L}_{\rho(\eta)}\SP{l, \xi} \ \ \text{(cf.\ equation\ \ref{def-xi})}\\
= \SP{l, \frac{d\eta}{dt}}+ \mathcal{L}_{\rho(\eta)}\SP{l, \eta}- d_{A}(l)(t, \gamma(t))\\
= \frac{d}{dt} \SP{ l(\gamma(t)), \eta(t, \gamma(t))}-
d_{A}(l)(a(t), \eta(t, \gamma(t))),
\end{split}
\]
and this clearly implies the desired formula. Now, (\ref{nice}) and (\ref{e1-form}) imply $i_{X_{\eta}}(df_{l}- \tilde{\sigma})= 0$. 
Using Cartan's formula $\mathcal{L}_X= di_X+ i_Xd$, we also deduce that
$\mathcal{L}_{X_{\eta}}(df_{l}- \tilde{\sigma})= -i_{X_{\eta}}(d\tilde{\sigma})$, 
which, as we mentioned above, is zero ($d\tilde{\sigma}$ was basic
with respect to the action of $\mathfrak{g}$). We deduce that
$(df_{l}- \tilde{\sigma})$ is a basic 1-form on $P(A)$, hence it
descends to a 1-form on $G(A)$, call it $\theta$. That $\theta$ is
multiplicative follows immediately from the fact that the
multiplication of $G(A)$ is defined by concatenation,
$(df_{l}- \tilde{\sigma})$ are defined by integrals, and
integration is an additive operation. Since the differential of $(df_{l}- \tilde{\sigma})$ is precisely the 2-form $- d\tilde{\sigma}$ that induces
$\omega$, we deduce that $d\theta= \omega$. Finally, for the uniqueness of the multiplicative 1-form with $d\theta= \omega$ and $\theta|_{A}= l$ we have to show that, if $\theta$ is closed an multiplicative and $\theta|_{A}= 0$, then $\theta= 0$. But this works exactly as for 2-forms (the argument being even simpler) for which we refer to Corollary 3.4 in \cite{CH}. 
\end{proof}

\subsection{The case of $G_{\omega}$}
\label{The case of G-omega}
We now derive the consequence that we need for prequantization. Let $\omega$ be a closed multiplicative $2$-form on
the Lie groupoid $G$, and we define
\[ \mathcal{P}_x(\omega)= Per(\omega|_{s^{-1}(x)}) , \mathcal{S}(\omega)= \mathbb{R}_{M}/\mathcal{P}(\omega).\]
Note that, if $c$ is the algebroid 2-cocycle induced by $\omega$ (see subsection \ref{From groupoids to algebroids}), then
$\mathcal{P}(\omega)$ and $\mathcal{S}(\omega)$ coincide with the similar bundles of $c$ (see \ref{The case of 2-cocycles}).
We also consider the induced algebroid $A_c$ (see subsection
\ref{The infinitesimal analogue}), and we put
\[ G_{\omega}:= G(A_c) .\]

\begin{corollary} 
\label{cor-2-form} Let $\omega$ be a closed multiplicative 2-form on an $s$-simply connected Lie groupoid,  
and assume that $G_{\omega}$ (or, equivalently, $\mathcal{P}(\omega)$, or $\mathcal{S}(\omega)$, cf. Theorem 
\ref{2-coc-int}) are smooth.

Then there is an extension of Lie groupoids
\[ \mathcal{S}(\omega)\rmap G_{\omega}\stackrel{\pi}{\rmap} G ,\]
and $G_{\omega}$ comes equipped with an invariant, 
multiplicative 1-form $\theta$ which is uniquely determined by the
equations
\[ d\theta= \pi^*\omega,\ \theta(\frac{d}{dt})= 1 .\]
\end{corollary}

\begin{proof} We use Theorem \ref{theorem-new} applied to the canonical cocycle $l_c$ on $A_c$ which transgresses $\pi^*(c)$
(see subsection \ref{The infinitesimal analogue}) to produce $\theta$. The last formula can be checked for the 1-form $df_{l}- \tilde{\sigma}$
on $P(A_c)$ which was used to construct $\theta$:  there is no contribution from $\tilde{\sigma}$ since it does not involve $t$,
while $df_{l}(\frac{d}{dt})= \mathcal{L}_{\frac{d}{dt}}(f_{l})$ clearly equals $1$. From the last two equations it also follows that
$\mathcal{L}_{\frac{d}{dt}}(\theta)= 0$, i.e. $\theta$ is invariant. 
\end{proof}

\subsection{Other consequences} 
Let us point out some immediate consequences of Theorem \ref{theorem-new}.
First of all, combining with Theorem \ref{th-CH}, we deduce:

\begin{corollary} Given an $s$-simply groupoid, there is a
bijection between multiplicative 1-forms on $G$ and pairs
$(\rho^*, l)$ where $\rho^*: A\rmap T^*M$ and $l\in C^1(A)$
satisfying:
\begin{eqnarray*}
(d_{A}l)(\alpha, \beta)& = & \SP{\rho^*(\alpha), \rho(\beta)}\\
\rho^{*}([\alpha, \beta]) & = &
\mathcal{L}_{\alpha}(\rho^{*}(\beta)) -
\mathcal{L}_{\beta}(\rho^{*}(\alpha))+
d\SP{\rho^*(\alpha),\rho(\beta)}
\end{eqnarray*}
\end{corollary}

Next, there is a version of our discussion in the presence of a ``3-form background'', i.e. a closed 3-form $\phi\in \Omega^3(M)$.
Given $\phi$, a form $\omega\in \Omega^2(G)$ is called $\phi$-closed if $d\omega= s^{*}\phi- t^{*}\phi$, and Theorem \ref{th-CH}
has a version that applies to $\phi$-closed multiplicative 2-forms \cite{CH}. Given such a pair $(\omega, \phi)$, we say it is multiplicatively exact if there is a pair $(\theta, \chi)$ where $\theta\in \Omega^1(G)$ is multiplicative and $\chi\in \Omega^2(M)$, such that
\begin{equation}
\label{eq-twisted-f} 
\omega= d\theta+ s^*\chi- t^*\chi,\ \phi= d\chi .
\end{equation}
Given $(\omega, \phi)$, an interesting problem is to make it multiplicatively exact by pulling it back along various groupoids maps. The prequantization problem (and its twisted versions \cite{LX}) is a particular case, and other examples come from Gauge equivalences of Dirac structures and moment maps (see \cite{CH} and the references therein). To state the result in this case, remark that $\chi\in \Omega^{3}(M)$ induces $\rho^*(\chi)\in C^{3}(A)$ by composition with $\rho$.

\begin{corollary} Given $\phi\in \Omega^3(M)$ closed and $\omega$ a $\phi$-closed multiplicative $2$-form on the $s$-simply connected Lie groupoid $G$, there is a bijection between
\begin{enumerate}
\item[(i)] pairs $(\theta, \chi)$ satisfying (\ref{eq-twisted-f}), with $\theta\in \Omega^1(G)$ multiplicative.
\item[(ii)] pairs $(l, \chi)$ where $l\in C^1(A)$ and $\chi\in \Omega^2(M)$ satisfy
\[ d_{A}(l)= c_{\omega}+ \rho^*(\chi), \ d\chi= \phi .\]
\end{enumerate}
\end{corollary}

\begin{proof} Apply theorem \ref{theorem-new} to $\omega- s^*\chi+ t^*\chi$. 
\end{proof}

\begin{remark}\rm \ 
\label{remark-Br}
(slightly speculative)  
Note the similarity between the formulas appearing in our discussions of multiplicative forms and some of those appearing in (several  different places of) Brylinski's book \cite{Br}. 
For instance, the 2-form on $P(A)$ induced by $\phi$ and needed in the reconstruction of multiplicative forms in the presence of a ``3-form background'' (denoted by $\omega_{\phi}$ in Section 5 of \cite{CH}, and first considered in \cite{CaXu}), when applied to $A= TM$, gives the 2-form $\beta_{\gamma}$ on the loop space described in \cite{Br} (see formulas (3-8) and (6-8) there). As $A= TM$ in this case, the groupoid in the discussion is just the pair groupoid of $M$. However, there are several other different groupoids that give rise to such striking similarities. For instance, the formulas in \cite{Br} for the canonical gerbe of a compact simple Lie group $G$ are very similar
to the ones appearing in the discussion of multiplicative forms on the {\it action groupoid} of $G$ with respect to the adjoint action \cite{CH}, which come from the Lie group-valued moment maps \cite{AMM} (this is closely related also to the results announced in \cite{BX}, and to work of P. Xu on Morita equivalences). Another case is that of groupoids induced by submersions that we describe next. 
First of all, any submersion $f: Y\rmap X$ induces a groupoid $Y\times_{X}Y$ over $Y$, similar to the pair groupoid (the source and target map are the projections). On the other hand, the submersion induces a filtration on $\Omega(Y)$, and we are interested here only in the spaces
\[ F_2\Omega^n(Y)= \{ \omega\in \Omega^2(M): i_{v_1}\ldots i_{v_{n-1}}\omega= 0\ \text{if\ the\ $v_i$'s\ are\ tangent\ to\ the\ fibers\ of\ $f$}\} .\]
We say that two forms $\omega$ and $\omega'$ are $F_2$-equivalent if their difference is in $F_2$. Assume that the fibers of $f$ are simply connected and connected. Then Theorem \ref{th-CH} applied to $Y\times_{X}Y$ immediately implies

\begin{corollary} There is a 1-1 correspondence between closed multiplicative 2-forms on $Y\times_{X}Y$ and $F_2$-equivalence classes of 2-forms forms $\beta\in \Omega^{2}(Y)$ with the property that $d\beta\in F^2$. 
\end{corollary}

Note that the condition appearing in the corollary is exactly the one appearing in the construction of connective structures in \cite{Br} (see Lemma 5.3.3). Even more, the data used in \cite{Br} for the construction of the gerbe under discussion (the closed relative 2-form $\omega$ there) is closely related to the space of closed multiplicative 2-forms on $Y\times_{X}Y$ modulo the space of
multiplicatively exact 2-forms (this follows from Theorem \ref{theorem-new}). This last ``coincidence'' seems to be related to the  ``bundle-gerbe picture'', in a more general context of ``bundle gerbes over groupoids''. 
\end{remark}

\section{Prequantization of groupoids}
\label{Prequantization of groupoids}

Throughout this section, $G$ is an $s$-simply connected Lie groupoid, and $\omega\in \Omega^2(G)$ is closed and multiplicative. 
The analogue of the group of periods from the classical situation is the {\it the
bundle of periods} of $\omega$, $\mathcal{P}(\omega)$ (see \ref{The case of G-omega}):
the fiber at $x$ is the group of periods of the restriction of $\omega$ to the $s$-fiber above $x$.
Recall also that $G_{\omega}$ stands for the monodromy group $G(A_{c})$ 
of the algebroid $A_c$ 
associated to the algebroid 2-cocycle $c$ which is the restriction of $\omega$ to $A$.
First of all, we have:

\begin{theorem}
\label{main-theorem}
The following are equivalent:
\begin{enumerate}
\item[(i)] $(G, \omega)$ is prequantizable.
\item[(ii)] $\mathcal{P}(\omega)\subset \mathbb{Z}_{M}$.
\end{enumerate}
Moreover, in this case $G_{\omega}$ and $\mathcal{S}(\omega)$ are smooth, there 
is an extension of Lie groupoids
\[ \mathcal{S}(\omega)\rmap G_{\omega} \stackrel{\pi}{\rmap} G ,\]
and $G_{\omega}$ comes equipped with an invariant multiplicative $1$-form $\theta$ satisfying $\pi^{*}(\theta)= d\omega$ and 
$\theta(\frac{d}{dt})= 1$.

Finally, the prequantization of $(\tilde{G}, \tilde{\theta})$ of $(G, \omega)$ will be unique, and it is the quotient of
$(G_{\omega}, \theta)$ by the action of $\mathbb{Z}_{M}/ \mathcal{P}(\omega)$.
\end{theorem}

As we mentioned at the beginning of the introduction, although the base manifold and the $s$-fibers of Lie groupoids 
are allays assumed to be Hausdorff, we do allow our groupoids to be non-Hausdorff. In the Hausdorff case however, we can slightly improve our conclusions. 

\begin{corollary}
\label{main-theoremH} 
Given an $s$-simply connected Hausdorff Lie groupoid $G$ and a closed multiplicative 2-form $\omega\in \Omega^2(G)$, the following are equivalent:
\begin{enumerate}
\item[(i)] $(G, \omega)$ is prequantizable.
\item[(ii)] $\mathcal{P}(\omega)\subset \mathbb{Z}_{M}$.
\item[(iii)] $\mathcal{P}(\omega)= M\times (k\mathbb{Z})$ for some integer $k$.
\item[(iv)] $\omega$ is integral.
\item[(v)] $Per(\omega)\subset \mathbb{Z}$.
\end{enumerate}
Moreover, in this case there is 
a central extension of (Hausdorff) Lie groupoids 
\[ \mathbb{R}/k\mathbb{Z}\rmap G_{\omega}\stackrel{\pi}{\rmap} G ,\]
and $G_{\omega}$ comes equipped with an invariant multiplicative $1$-form $\theta$ satisfying $\pi^{*}(\theta)= d\omega$ and 
$\theta(\frac{d}{dt})= 1$. 

Finally, the prequantization of $(G, \omega)$ will be unique, it is Hausdorff, and it is obtained
from $(G_{\omega}, \theta)$ by dividing out by the action of the cyclic group of order $k$, $\mu_{k-1}\subset S^1$.
\end{corollary}

We also point out the following immediate consequence of the corollary 

\begin{corollary} A closed multiplicative 2-form on a Hausdorff groupoid $G$ is exact if and only if its restriction
to each $s$-fiber is exact.
\end{corollary}

\begin{proof} (of Theorem \ref{main-theorem}) 
The $s$-fibers of a prequantization of $(G, \omega)$ are classical prequantizations of the $s$-fibers of $G$, hence
condition (i) in the theorem implies (ii). Assume now that (ii) holds. First of all, using Theorem \ref{2-coc-int},
$A_c$ is integrable, hence using also Theorem \ref{cor-2-form}, we obtain the Lie groupoid extension from the statement
and the 1-form $\theta$. 
Next, since $G_{\omega}$ is principal with respect to the action of $\mathcal{S}(\omega)$,
and since $\mathbb{Z}_{M}/ \mathcal{P}(\omega)$ is the kernel of the canonical map $\mathcal{S}(\omega)\rmap S^{1}_{M}$,
it follows that $\mathbb{Z}_{M}/ \mathcal{P}(\omega)$ is a smooth bundle of discrete Lie groups and the 
quotient $\tilde{G}_{\omega}$ of $G_{\omega}$ with respect to the action of this bundle will be a Lie groupoid (of the same dimension as $G_\omega$). The invariance of
$\theta$ implies that it descends to a 1-form $\tilde{\theta}_{\omega}$ on $\tilde{G}_{\omega}$ with the same properties as $\theta$. 
Hence $(\tilde{G}_{\omega}, \tilde{\theta}_{\omega})$ is a prequantization of $(G, \omega)$. 
If $(\tilde{G}, \tilde{\theta})$ is another prequantization, we already know that
$\tilde{G}$ will integrate $A_c$ (see Lemma \ref{first-step}), and from this we obtain a map 
$f: G_{\omega}= G(A_c)\rmap \tilde{G}$.
This follows from the general properties of monodromy groupoids \cite{CrFe1} and it can be easily described
using the correspondence between $A_c$-paths and paths in $\tilde{G}$ (the $D^R$ induced by the formula (\ref{bij-paths})):
 given $[a]\in G(A_c)$, we choose $g:I\rmap \tilde{G}$ sitting in an $s$-fiber such that $D^{R}(g)= a$, and $f([a])= g(1)$. 
Since this map clearly commutes with  the projections into $G$, both $\theta$ and $f^{*}\tilde{\theta}$ are multiplicative
and transgress $\pi^*\omega$, using also Theorem \ref{theorem-new}, we deduce that $f^{*}\tilde{\theta}= \theta$.
Since $\theta(\frac{d}{dt})= 1$ and $\tilde{\theta}(\frac{d}{dt})= 1$ we deduce that the restriction of $f$ to $\mathcal{S}(\omega)$,
$f: \mathcal{S}(\omega)\rmap S^{1}_{M}$ is the canonical map that sends $[r]\in \mathbb{R}/\mathcal{P}_x(\omega)$ into
the class of $r$ in $\mathbb{R}/\mathbb{Z}$. This shows that the kernel of $f$ is indeed $\mathbb{Z}/\mathcal{P}(\omega)$ hence
$f$ descends to give an isomorphism $\tilde{G}_{\omega}\rmap \tilde{G}$ (its inverse is smooth because $f$ is a submersion).
\end{proof}

\begin{proof} (of Corollary \ref{main-theoremH}) 
We already know that  (i) and (ii) are equivalent,
and clearly (iv) implies (v), (v) implies (ii), and (iii) implies (iii). Also, since a prequantization of the groupoid $G$ is in particular a prequantization in the   
classical sense, (i) implies (iv) at least if we know that the prequantization is Hausdorff. But, in general, we have the following immediate remark

\begin{lemma} Given an extension $\mathcal{S}\rmap \tilde{G}\rmap G$ of Lie groupoids with $G$ Hausdorff, $\tilde{G}$ is Hausdorff if and only if $\mathcal{S}$. In particular, prequantizations of a Hausdorff groupoid are always Hausdorff.
\end{lemma}

Hence we are left with proving that (ii) implies (iii). This follows immediately from the following property 

\begin{lemma} If $\omega$ is a closed multiplicative 2-form on a Hausdorff Lie groupoid $G$, then $\mathcal{P}(\omega)$ has enough smooth sections, i.e., for any $x\in M$ and any $s_0\in \mathcal{P}_{x}(\omega)$, there is a smooth local section of $\mathbb{R}_{M}$ around $x$ such that $s(x)= s_0$ and $s(y)\in \mathcal{P}_y(\omega)$ for all $y$.
\end{lemma}

For the proof, one writes $s_0= \int_{g} \omega|_{s^{-1}(x)}$ for some $g: S^2\rmap s^{-1}(x)$ representing a 2-homotopy class. Since
$s: G\rmap M$ is a submersion between Hausdorff manifolds, $S^2$ is simply connected and compact, it follows that we can deform $g$ to a smooth family of spheres $g_y$ in $s^{-1}(y)$ for $y$ close enough to $x$, such that $g_x= g$. Then
$s(y)= \int_{g_y}\omega|_{s^{-1}(y)}$ is a smooth local section of $\mathbb{R}_{M}$ passing through $s_0$.  
\end{proof}

\begin{footnotesize}

\end{footnotesize}

\end{document}